\documentclass[reqno,10pt]{article}
\usepackage{hyperref}
\usepackage{amsmath}
\usepackage{amssymb}
\usepackage{mathrsfs}

\newtheorem{thm}{Theorem}[section]
\newtheorem{lemma}[thm]{Lemma}

\newtheorem{corol}[thm]{Corollary}
\newtheorem{propos}[thm]{Proposition}
\newtheorem{rema}{Remark}[section]
\def\bp{\begin{propos}}
\def\ep{\end{propos}}
\def\bt{\begin{thm}}
\def\et{\end{thm}}
\def\bco{\begin{corol}}
\def\eco{\end{corol}}
\def\bl{\begin{lemma}}
\def\el{\end{lemma}}
\def\br{\begin{rema}}
\def\er{\end{rema}}
\def\be{\begin{equation}}
\def\ee{\end{equation}}
\def\ba{\begin{array}}
\def\ea{\end{array}}
\def\bena{\begin{eqnarray}}
\def\eena{\end{eqnarray}}

\def\P{{\mathbb P}}
\def\E{{\mathbb E}}
\def\R{{\mathbb R}}

\def\V{{\rm Vol}}
\def\1{I}
\def\imath{\textbf{i}}
\def\jmath{\textbf{j}}
\def\chi{\zeta}

\def\fD{{\cal D}}
\def\hE{{\cal E}}
\def\fE{{e}}

\def\hB{{\mathscr B}}

\def\fG{{\mathscr G}}

\def\fL{{\cal L}}
\def\fP{{\mathscr P}}

\def\a{{\alpha}}

\def\b{\binom}

\def\QED{\hfill$\square$\vskip 3mm}

\def\Dp{\displaystyle}
\def\Df{\Dp\frac}

\def\({\left(}
\def\){\right)}
\def\ln{\lg}

\begin{document}

\title{Mixing Time of Random Walk on Poisson\\ Geometry Small World
\\[5mm]
\footnotetext{AMS classification (2000): Primary 05C 80; secondary 60K 35£¬ 60K 37 } 
\footnotetext{Research supported in part by the Natural Science Foundation of China (under
grants 11271356, 11471222) and the Foundation
of Beijing Education Bureau (under grant KM201510028002) }}

\author{Xian-Yuan Wu}\vskip 10mm
\date{}
\maketitle

{\begin{center} \begin{minipage}{10cm}
{\bf Abstract}: This paper focuses on the problem of modeling for {\it small world effect} on complex networks. Let's consider the supercritical Poisson continuous percolation on $d$-dimensional torus $T^d_n$ with volume $n^d$.
By adding ``long edges (short cuts)" randomly to the largest percolation cluster, we obtain a random graph $\fG_n$.
In the present paper, we first prove that the diameter of $\fG_n$ grows at most polynomially fast in $\ln n$ and we call it the
{\it Poisson Geometry Small World}. Secondly, we prove that the random walk on $\fG_n$ possesses the {\it rapid mixing} property, namely, the random walk mixes in time at most polynomially large in $\ln n$.

\end{minipage}
\end{center}}

 \vskip 5mm
\section{Introduction and statement of the results}
\renewcommand{\theequation}{1.\arabic{equation}}
\setcounter{equation}{0}
{\it Small world effect}, the fact that the diameters of most networks are considerably smaller than their sizes,
 is one of the most important features of real-world complex networks. The existence of small world effect had been speculated upon in a remarkable
 short story by Karinthy \cite{K} in 1929. In 1960s, Milgram \cite{M,TM} carried out his famous ``small-world" experiments, in which letters passed
 from person to person were able to reach a designated target individual within six steps, and which finally led to the popular concept of
 the ``six degrees of separation" \cite{G}. Recent influential studies on small world effect perhaps started with the work of Watts and Strogatz
 published in 1998 \cite{WS}. From then on, people were much more interested in studying the structure features (including {\it small world effect},
 {\it scale-free property} and {\it navigability}, etc.) of complex networks.
Nowadays, the small world effect has been studied and verified directly in a large number of different networks, see \cite[Table 3.1]{N} and the references therein.

What are the underlying causes which make most networks small worlds? To answer this question, many models have been introduced and studied by physicists and mathematicians. The most important ones include the Bollob\'as and Chung small world model \cite{BC} (BC small world), the Newman and Watts small world model \cite{NW} (NW small world) and the Watts and Strogatz small world model \cite{WS} (WS small world). Actually, all these models were introduced to reveal such a fact that adding ``long edges (short cuts)" to a regularly constructed (lattice-like) graph will make the resulted graph a small world, and we will call it the {\it adding-long-edges} mechanism. It should be noted that in the above three models, only \cite{BC} provided rigorous mathematical results. For other mathematical results on small world effect, one may refer to \cite{BR,CL}. In \cite{BR}, an evolving random graph process, which is called the `LCD' model, was introduced to model the evolution of real-world complex networks. It was proved that, while the model ultimately possess a power law degree distribution, the model also exhibit small world effect. It seems that a mechanism other than the one working in \cite{BC,NW} and \cite{WS} makes the `LCD' model a small world.

The present paper will introduce a new model to study the small world effect of real-world complex networks. Precisely speaking, a new and more appropriate model will be introduced to explain the {\it adding-long-edges} mechanism mentioned in the past paragraph, and we will call it the {\it Poisson Geometry Small World}.

First of all, let's recall the Poisson continuous percolation on $\R^d$, $d\geq 2$. Let $\fP$ denote the homogeneous Poisson process of rate $1$ on $\R^d$. Given $r>0$, define the ($r$-){\it clusters} on $\fP$ to be the connected components of the union of the balls of radius $r$ centered at the points of $\fP$, here radius $r$ is relative to the usual Euclidean metric. We call {\it percolation occurs} when an unbounded cluster exists. Define the percolation probabilities $\theta(r)$ and $\tilde\theta(r)$ as follows: let $\theta(r)$ denote the probability that there is an unbounded cluster on $\fP$ containing the origin $0$, and let $\tilde\theta(r)$ denote the probability that there is an unbounded cluster that intersects the ball of radius $r$ centered at $0$. Then $\theta(\cdot)$ and $\tilde\theta(\cdot)$ are non-decreasing and for each $r$, $\theta(r)$ and $\tilde\theta(r)$ are either both zero or both strictly positive. Define the critical value $r_c=r_c(d)$ by $r_c:=\inf\{r:\theta(r)>0\}$. Assume that $d\geq 2$ and in this case $0<r_c<\infty$. For more details on this model, one may refer to \cite{Gr,MR,P}.


Let $B(n)$ denote the cube $[0,n]^d$, and set $\fP_n:=\fP\cap B(n)$, a Poisson process of rate $1$ on $B(n)$. Given $r>0$, define the {\it clusters} on $\fP_n$ to be the connected components of the {\it intersection} of $B(n)$ and the union of the balls of radius $r$ centered at the points of $\fP_n$. A cluster $C$ on $\fP_n$ is called {\it crossing} for $B(n)$, if $C$ intersects all of the $2d$ faces of $B(n)$.

Let $T^d_n$ denote the $d$-dimensional torus obtained from $B(n)$ by cohering its opposite faces, and let $\fP^T_n$ denote the Poisson process of rate $1$ on $T^d_n$. Given $r>0$, define the {\it clusters} on $\fP^T_n$ to be the connected components of the union of the balls of radius $r$ centered at the points of $\fP^T_n$. Here radius $r$ is relative to the metric on $T^d_n$, the metric naturally inherited from the Euclidean metric on $\R^d$. We shall be interested in the clusters on $\fP^T_n$ for $r>r_c$.

For any measurable $A\subset T^d_n$ or $B(n)$, let $|A|$ denote its cardinality if $A$ is finite or countable, or its volume (Lebesgue measure) otherwise. For any measurable $A\subset B(n)$ or $\R^d$, let diam$(A)$ denote its diameter, i.e. diam$(A)=\sup\{||x-y||:x,y\in A\}$, with $||\cdot||$ the Euclidean norm.
A cluster $C$ on $\fP_n$ or $\fP^T_n$ is called the {\it largest cluster} if $|C|$ reaches the maximum. By the property of Poisson process, for any $r>0$, there is asymptotically almost surely as $n\rightarrow\infty$ a {\it unique} largest cluster on $\fP_n$ and $\fP^T_n$. Let $C_{\rm max}$ denote the unique largest cluster on $\fP^T_n$. To any cluster $C$ on $\fP^T_n$, we associate a graph $G_C=(V_C,E_C)$, where $V_C$ is the vertex set defined by $V_C:=C\cap \fP^T_n$ and $E_C$ is the set of edges which connect all $V_C$ vertex pairs lying in distance $2r$ from each other. Let $G_n=(V_n,E_n)$ denote the graph associated to the largest cluster $C_{\rm max}$.

Let $d^T_\infty(\cdot,\cdot)$ denote the $l_\infty$ metric on $T^d_n$ inherited from the usual $l_\infty$ metric $d_\infty(\cdot,\cdot)$ on $\R^d$ defined by $d_\infty(x,y):=\Dp\max_{1\leq s\leq d}|x_s-y_s|$ for any $x,y\in\R^d$. For any given constants $\a,\beta,\sigma$ and $\chi$ satisfying $0<\a<\beta <1/2$, $\sigma>0$ and $\chi\in \R$, we define a random graph $\fG_n=\fG_n(\a,\beta,\sigma,\chi)$ from $G_n$ as follows: for any $u,v\in V_n$, if $\a n\leq d^T_\infty(u,v)\leq \beta n$, then we connect $u$ and $v$ independently by a ``long edge" with probability \be p_n=\sigma n^{-d}\ln^{\chi}n;\ee otherwise, we do nothing. Let $\hE_n$ denote the new edge set with long edges, and let $\fG_n=(V_n,\hE_n)$.

Now, we have finished the definition of the model. We hope $\fG_n$ can be qualified to model the small world effect of some kind of real-world complex networks. In fact, $\fG_n$ can be seen as a {\it higher dimensional} and {\it random based} version of the NW small world proposed in \cite{NW}. Recall that the NW small world started from a {\it ring lattice} with $n$ vertices, then a Poisson number of shortcuts (i.e. long edges) with mean $\Theta(n)$ are added and attached to randomly chosen pairs of sites. In our model $\fG_n$, firstly, the random graph $G_n$ plays the same role as the {\it ring lattice} played in the construction of the NW small world; secondly, we add a Binomial number of long edges with mean $\Theta(n^{2d}p_n)$ and we add them in such a way that no {\it double edge} appears in $\fG_n$; finally, in $\fG_n$, only links between two vertices at Euclidean distance $\Theta(n)$ are treated as ``long edges".

In the past paragraph we have used the notation $\Theta(b_n)$, in fact, we use $a_n=\Theta(b_n)$ to denote $cb_n\leq a_n\leq C b_n$ for some $0<c<C<\infty$. For convenience, in this paper we also use $a_n=O(b_n)$ to denote $a_n\leq Cb_n$ for some $C>0$ and use $a_n=\Omega(b_n)$ to denote $a_n\geq cb_n$for some $c>0$.

We will first study the diameter of $\fG_n$. Recall that in a graph $G$, the distance $\fD_G(u,v)$ between two vertices $u$ and $v$ is the length (number of edges) of the shortest path between them, and the {\it diameter} diam$(G)$ of a connected graph $G$ is the maximum distance between two vertices.

\bt\label{th1}
Suppose $r>r_c$. Then
\begin{description}
  \item[] (i)  for any $0<\a<\beta<1/2$, $\chi\leq 1/(1-d)$, and for $\sigma>0$ small enough, there exists constant $C>0$ such that
  \be\label{3}\lim_{n\rightarrow\infty}\P({\rm diam}(\fG_n)\geq C\ln^{(1-\chi)/d} n)=1;\ee
  \item[] (ii) for any $0<\a<\beta<1/2$ with $(2\beta)^d-(2\a)^d>1/2$, for $\chi> 1$ and $\sigma>0$, or for $\chi=1$ and $\sigma$ is large enough, there exists constant $C_1>0$ such that
\be\label{1}\lim_{n\rightarrow\infty}\P({\rm diam}(\fG_n)\leq C_1\ln^{2} n)=1.\ee
\end{description}
\et

\br It seems that our setting on ``long edge" is {REASONABLE}! Obviously, if only shorter edges, for example with length $n^{1-\epsilon}$, are added, then the diameter of the resulted graph grows at least fast as $n^\epsilon$, and the resulted graph does not exhibit the small world effect. On the other hand, Theorem~\ref{th1} indicates that, to make the resulted graph a small world, adding such shorter edges is not necessary.\er

\br\label{r} In the definition of $\fG_n$, we only used Poisson points in the largest percolation cluster $C_{\rm max}$ as its node set. In fact,
 we may obtain a more complicated random graph $\bar\fG_n$ by connecting each $G_C$, $C\not=C_{\rm max}$, to $\fG_n$ with a additional shortest edge (in Euclidean distance) between $V_C$ and $V_n$.
 By Theorem~\ref{th1} and Proposition~\ref{p2} below, $\bar \fG_n$ is also a small world. \er

\br The problem for giving upper or lower bound to ${\rm diam}(\fG_n)$ when $1/(1-d)< \chi< 1$ remains open. It is a pity that we can not give both lower and upper bounds to ${\rm diam}(\fG_n)$ for any given $\chi$. Furthermore, it seems that the bounds given by Theorem~\ref{th1} are sub-optimal. Note that, in \cite{ZW}, for a modified NW model (where the d-dimensional lattice torus takes the place of $G_n$), its diameter is bounded from below and above by power functions of $\ln n$ when $0\leq\chi\leq 1$. \er

As Newman noted in \cite{N}, the ultimate goal of the study of the structure of networks is to understand and explain the workings of systems built upon those networks. Clearly, random walks on networks are just the simplest (but important) workings of systems built upon networks. At the present paper, we will next study the {\it mixing time} of random walk on $\fG_n$. In probability theory, the mixing time of a Markov chain is the time until the Markov chain is ``close" to its steady state distribution. The concept of mixing times was presented to a wider-range
audience by Aldous and Diaconis in 1986 \cite{AD}. Since then, both the mathematical theory and its interactions with computer science and statistical physics have been developed tremendously. While mathematical theory mainly focuses on how mixing times change as a function of the size of the structure underlying the chain \cite{AF,LPW}, the most developed theory manages randomized algorithms for NP-Complete algorithmic counting problems in computer science, see \cite{S} {\it etc}. For basic concepts on mixing time and related problems on mathematics and statistical physics, one may refer to \cite{LPW} and the references therein. For mixing time of random walk on complex networks, one may refer to \cite{D,LPW}.

In a graph $G=(V,E)$, for any $u,v\in V$, let $d_G(u)$ be the the degree of $u$ in $G$, and write $u\sim v$ if $u$ and $v$ are neighbors in $G$. Let $\Delta(G)$ denote the {\it maximum degree} of $G$, i.e. $\Delta(G):=\max\{d_G(u):u\in V\}$.

For any $u,v\in V_n$, we define a transition kernel by $P(u,u)=1/2$, $P(u,v)=1/2d_{\fG_n}(u)$ if $u\sim v$ and $P(u,v)=0$ otherwise. A discrete time Markov chain $\{X_t:t\geq 0\}$ on $V_n$ with transition kernel $(P(u,v))$ is called the {\it lazy random walk} on $\fG_n$. Note that $\pi(u):=d_{\fG_n}(u)/D$ where $D=\sum_{v\in V_n}d_{\fG_n}(v)$, defines a reversible stationary distribution of $\{X_t\}$ since $\pi(u)P(u,v)=1/2D=\pi(v)P(v,u)$. By the basic theory of Markov chains, for any initial state $u\in V_n$, the distribution of $X_t$, i.e. $P^t(u,\cdot):=\P(X_t\in\cdot\mid X_0=u)$, converges weakly to $\pi$ as $t\rightarrow\infty$. To measure convergence to equilibrium, we will use the {\it total variation distance} $$||P^t(u,\cdot)-\pi||_{TV}:=\frac 12\sum_{v\in V_n}|p^t(u,v)-\pi(v)|.$$ The mixing time of $\{X_t:t\geq 0\}$ is defined by \be\label{2} T_{\rm mix}:=\min\left\{t:\max_{u\in V_n}||P^t(u,\cdot)-\pi||_{TV}<1/e\right\}.\ee

The second result of the present paper is about $T_{\rm mix}$ and we state it as follows

\bt\label{th2}
Suppose $r>r_c$. Then
\begin{description}
  \item[] (i)  for any $0<\a<\beta<1/2$, $\chi\leq 1/(1-d)$ and for $\sigma>0$ small enough, there exists constant $C_2>0$ such that
  \be\label{3}\lim_{n\rightarrow\infty}\P(T_{\rm mix}\geq C_2\ln^{(1-\chi)/d} n)=1;\ee
  \item[] (ii) for any $0<\a<\beta<1/2$ with $(2\beta)^d-(2\a)^d>1/2$, for $\chi> 1$ and $\sigma>0$, or for $\chi=1$ and $\sigma$ is large enough, there exists constant $C_3>0$ such that
\be\label{4}\lim_{n\rightarrow\infty}\P(T_{\rm mix}\leq C_3\ln^3 n)=1.\ee
\end{description}
\et

\br Theorem~\ref{th2} provides a good random sampling method to get $\pi$, the stationary distribution of $X_t$, which exhibits important structural properties of $\fG_n$. Because $\bar\fG_n$ is constructed from $\fG_n$ in a simple way as given in Remark~\ref{r}, one may also understand $\bar\fG_n$ through the stationary distribution $\pi$.\er

The rest of the paper is arranged as follows. In Section 2, we study the geometry of the largest cluster of supercritical Poisson continuous percolation on $T^d_n$. In Section 3, we bound the maximum degree of $\fG_n$ from above, and bound the {\it Cheeger constant} of $\fG_n$ and the {\it conductance} of random walk on $\fG_n$ from below. Finally, we prove Theorems~\ref{th1},~\ref{th2} in Section 4.

\vskip 5mm
\section{Geometry of the largest cluster}
\renewcommand{\theequation}{2.\arabic{equation}}
\setcounter{equation}{0}
In this section, based on the large deviation results of continuous percolation given by Penrose and Pisztora \cite{PP}, we characterize the geometry of the largest cluster on $\fP^T_n$.

Suppose $r>r_c$. Then there is almost surely a {\it unique} unbounded cluster $C_\infty$ on $\fP$ (see \cite{MR}). Therefore it is natural to expect that for large $n$, there is likely to be a {\it big} cluster on $\fP^T_n$ containing a proportion $\theta(r)$ of $T^d_n$. In fact, for supercritical continuous percolation on $B(n)$, the geometry of the largest cluster has been well studied in \cite[Theorem 1]{PP}. Here we state a subtly simplified version of the theorem, which is enough for our use, as the following proposition.
\bp\label{p1} Suppose $r>r_c$, and $0<\epsilon<1/2$. Let $E(n)$ be the event that (i) there is a unique cluster $C_b(B(n))$ on $\fP_n$ containing more than $\epsilon\tilde\theta(r)n^d$ points of $\fP_n$, (ii) \be(1-\epsilon)\theta(r)\leq n^{-d}|C_b(B(n))|\leq(1+\epsilon)\theta(r)\ee and
\be(1-\epsilon)\tilde\theta(r)\leq n^{-d}|C_b(B(n))\cap\fP_n|\leq(1+\epsilon)\tilde\theta(r),\ee (iii) $C_b(B(n))$ is crossing for $B(n)$, and (iv) $C_b(B(n))$
is part of the unbounded cluster $C_\infty$. Then, there exist $c_1>0$ and $n_0$ such that
\be\P(E(n))\geq 1-\exp(-c_1n^{d-1}), \ \ n\geq n_0.\ee \ep Note that items (i) and (ii) are translated to our fashion by using
the {\it scaling relation} of continuous percolation. The {\it scaling relation} of continuous percolation tells such a fact that, under
a $r'/r$-times magnifying glass, a system with parameter $(\lambda,r)$ is just the system with parameter $(\lambda',r')$ with
$\lambda'=\lambda(r/r')^d$. Where $\lambda$ is the rate of the Poisson process and $r$ is the radius of the concerned balls.


A more fundamental result for supercritical continuous percolation on $\fP_n$ was also studied in \cite[Proposition 2]{PP}. Through the scaling relation, we obtain the following proposition.

\bp\label{p2} Suppose $r>r_c$. Suppose $\{\phi_n:n\geq 1\}$ is increasing with $\phi_n/\ln n\rightarrow\infty$ as $n\rightarrow\infty$, and with $\phi_n<n$ for all $n$. Let $E'(n)$ be the event that (i) there is a unique cluster on $\fP_n$ that is crossing for $B(n)$, and (ii) no other cluster on $\fP_n$ has diameter greater than $\phi_n$. Then there exists a constant $c_2>0$ such that for all large enough $n$, \be\P(E'(n))\geq 1-\exp(-c_2\phi_n).\ee \ep

Based on Propositions~\ref{p1} and~\ref{p2}, we obtain the geometry of $C_{\rm max}$, the largest cluster on $\fP^T_n$ as stated in the following four lemmas.

\bl\label{l1} Suppose $r>r_c$, and $0<\epsilon<1/2$. Then there exists a constant $c_3>0$ such that for all large $n$
\be\label{11.1} \P\((1-\epsilon)\tilde\theta(r)\leq n^{-d}|V_n|\leq (1+\epsilon)\tilde\theta(r)\)\geq 1-\exp(-c_3n^{1/2}).\ee\el
{\it Proof.} Let $C_b(B(n))$ be the largest cluster on $\fP_n$. Let $C_1,C_2,\ldots,C_l$ denote all the other clusters which intersect the boundary of $B(n)$. Then by Proposition \ref{p1}, for any $\epsilon>0$, there exist $c_1>0$ and $n_0$ such that
\be\label{l1.01}
\ba{rl}&\P\(\ba{rl}1)&C_b(B(n))\subset C_{\rm max}\subset C_b(B(n))\cup\left\{\cup_{i=1}^lC_i\right\}, \ {\rm and}\\[1mm]
2)&(1- \epsilon/2)\tilde\theta(r)\leq n^{-d}|C_b(B(n))\cap \fP_n|\leq (1+\epsilon/2)\tilde\theta(r)\ea \)\\[6mm]
&\hskip3mm \geq 1-\exp(-c_1 n^{d-1}),\ \ \forall\ n\geq n_0.\ea\ee
On the other hand, by Proposition~\ref{p2}, let $\phi_n=n^{\frac 12}$, then for some $c_2>0$
\be\label{l1.02}
\P\(\cup_{i=1}^lC_i\subset B(n)\cap U_{n^{1/2}}(B(n)^c)\)\geq 1-\exp(-c_2n^{\frac 12}),\ee where $U_r(A)$ denotes the $r$ neighborhood of $A$, i.e. $U_r(A)=\{x:||x-y||\leq r\ {\rm for\ some}\ y\in A\}$. Let $N(n)$ denote the number of Poisson points in $B(n)\cap U_{n^{1/2}}(B(n)^c)$, a Poisson random variable with mean $\Theta(n^{d-1/2})$, then
\be\label{l1.04} |C_b(B(n))\cap \fP_n|\leq |V_n|\leq |C_b(B(n))\cap\fP_n|+N(n).\ee

By our large deviation result for Poisson distribution, if $Z\sim P(\mu)$, i.e. $Z$ obeys the Poisson distribution with mean $\mu$, then
\be\label{l1.05}\P(Z\geq z\mu)\leq \exp(-\gamma(z)\mu),\ {\rm where}\ \gamma(z)=z\ln z-z+1,\ z>1.\ee
So, there exists $c'>0$ such that
\be\label{l1.03}\P\(N(n)\geq n^{d-1/6}\)\leq \exp({-c'n^{d-1/6}}).\ee Combining (\ref{l1.01}), (\ref{l1.02}), (\ref{l1.04}), and (\ref{l1.03}), we have, for some $c_3>0$ and $n_1$
$$
\P\((1-\epsilon)\tilde\theta(r)\leq n^{-d}|V_n|\leq (1+\epsilon)\tilde\theta(r)\)\geq 1-\exp(-c_3n^{1/2})$$ for all $n\geq n_1$.\QED

\bl\label{l2}Suppose $r>r_c$, $0<\epsilon<1/2$ and $0<\a<\beta<1/2$. Then there exists a constant $c_4>0$ such that for all large $n$
\be\ba{ll}\label{l2.1}&\P\((1-\epsilon)\Gamma\leq n^{-d}|\Lambda_n(u)|\leq (1+\epsilon)\Gamma\ {\rm for\ all}\ u\in V_n \)\\[1mm]&\ \ \geq1-\exp(-c_4n^{1/2}),\ea\ee where $\Gamma=[(2\beta)^d-(2\a)^d]\tilde\theta(r)$ and $\Lambda_n(u)=\left\{v\in V_n:\a n\leq d^T_\infty(u,v)\leq\beta n\right\}$.
\el

{\it Proof.} For any $x\in T^d_n$, $\xi>0$, let $B^{T_n}_x(\xi)$ denote the box in $T^d_n$ centered at $x$ with side length $2\xi$. Let $A^\beta_\a(x,n)=B^{T_n}_x(\beta n)\setminus B^{T_n}_x(\a n)$. We will first prove that, under the assumption of the lemma, there exist $c_5>0$ and $n_2$ such that
\be\label{l2.01}\P\((1-\epsilon)\Gamma\leq n^{-d}|A^\beta_\a(x,n)\cap V_n|\leq(1+\epsilon)\Gamma\)\geq 1-\exp(-c_5n^{1/2})\ee for all $n\geq n_2$.

For any given $\epsilon'>0$ small enough, applying Proposition~\ref{p1} to the continuous percolation on boxes $B^{T_n}_x(\a n)$, $B^{T_n}_x(\beta n)$ and $B^{T_n}_x(\frac 12 n)$ respectively, we obtain
\be\label{l2.02}
\ba{rl}&\P\(\ba{rl}1)&C_b(B^{T_n}_x(\a n))\subset C_b(B^{T_n}_x(\beta n))\subset C_b(B^{T_n}_x(\frac 12 n))\subset C_{\rm max}, \ {\rm and}\\[1mm]
2)&(1-\epsilon')\tilde\theta(r)\leq(2\a n)^{-d}|C_b(B^{T_n}_x(\a n))\cap\fP^T_n|\leq (1+\epsilon')\tilde\theta(r), \ {\rm and}\\[1mm]
3)&(1-\epsilon')\tilde\theta(r)\leq(2\beta n)^{-d}|C_b(B^{T_n}_x(\beta n))\cap\fP^T_n|\leq (1+\epsilon')\tilde\theta(r)\ea \)\\[8mm]
&\hskip3mm \geq 1-\exp(-c n^{d-1})\ea\ee for some $c>0$ and all large $n$. Note that $\epsilon'$ is asked to be small to guarantee that the above item 1) hold by items (i) and (iv) of Proposition \ref{p1}.

By Proposition~\ref{p2} and the large deviation inequality (\ref{l1.05}), the number of all Poisson points in $[C_b(B^{T_n}_x(\beta n))\cap B^{T_n}_x(\a n)]\setminus C_b(B^{T_n}_x(\a n))$ and $[C_b(B^{T_n}_x(\frac 12 n))\cap B^{T_n}_x(\beta n)]\setminus C_b(B^{T_n}_x(\beta n))$ does not exceed $n^{d-1/6}$ with probability at least $1-\exp(-c'n^{1/2})$ for some $c'>0$ and all large $n$. This, together with (\ref{l2.02}), implies that
\be\label{l2.03}
\ba{rl}&\P\(\ba{rl}1)&(1-\epsilon')\tilde\theta(r)\leq (2\a n)^{-d}|B_x^{T_n}(\a n)\cap V_n|\leq (1+2\epsilon')\tilde\theta(r), \ {\rm and}\\[1mm]2)&(1-\epsilon')\tilde\theta(r)\leq (2\beta n)^{-d}|B_x^{T_n}(\beta n)\cap V_n|\leq (1+2\epsilon')\tilde\theta(r)\ea \)\\[6mm]
&\hskip3mm \geq 1-\exp(-c''n^{1/2})\ea\ee
for some $c''>0$ and all large $n$.

Take $\epsilon'$ small enough such that $\epsilon'\leq{\epsilon(\beta^d-\a^d)}/2(\beta^d+\a^d)$, (\ref{l2.01}) follows from (\ref{l2.03}).

Fix some integer $m$. For any ${\imath}=(i_1,i_2,\ldots,i_d)\in\{0,1,2,\ldots, m-1\}^d$, let $B^m_\imath(n)=\frac nm\imath+[0,\frac nm]^d$. Suppose that $(1-\epsilon)\Gamma\leq n^{-d}|\Lambda_n(u)|\leq (1+\epsilon)\Gamma$ does not hold for some $u\in V_n$. Then there exist $\imath\in\{1,2,\ldots,m-1\}^d$ and $x(=u)\in B^m_\imath(n)$ such that
\be\label{l2.05}(1-\epsilon)\Gamma\leq n^{-d}|\Lambda_n(x)|\leq (1+\epsilon)\Gamma\ {\rm does\ not\ hold}.\ee Denote by $x_\imath$ the center of the box $B^m_\imath(n)$, and rewrite $B^m_\imath(n)$ as $B^{T_n}_{x_\imath}(\frac n{2m})$.

Let's consider the boxes $B^{T_n}_{x_\imath}(\beta n\pm\frac n{2m})$ and $B^{T_n}_{x_\imath}(\a n\pm\frac n{2m})$. The difference between $|\Lambda_n(u)|$ and $|A^{\beta+1/2m}_{\a-1/2m}(x_\imath,n)\cap V_n|$
can not exceed $M=M(n,\a,\beta,m)$, the number of Poisson points in $A^{\beta+1/2m}_{\beta-1/2m}(x_\imath,n)$
and $A^{\a+1/2m}_{\a-1/2m}(x_\imath,n)$. Suppose $m_1$ is large enough and $a_1>0$ is small enough. Then, for $m\geq m_1$ and $0<a\leq a_1$, (\ref{l2.05}), together with $M\leq a n^d$, implies
\vskip-4mm
\be\label{l2.04}(1-\epsilon/2)\Gamma'\leq n^{-d}|A^{\beta'}_{\a'}(x_\imath,n)\cap V_n|\leq (1+\epsilon/2)\Gamma' \ {\rm does\ not\ hold,} \ee where $\Gamma'=\tilde\theta(r)\cdot[(2\beta')^d-(2\a')^d]$ and $\beta'=\beta+1/2m$, $\a'=\a-1/2m$.

Now, we take $m\geq m_1$ such that $\E(M)<{a_1}n^d/2 $. By the large deviation inequality (\ref{l1.05}), the probability that $M>a_1 n^{d}$ is less than $\exp(-\Omega(n^d))$. Then, by applying (\ref{l2.01}) to (\ref{l2.04}), we have
$$\P\(\ba{rl}&\hskip -5mm{\rm there \ exist }\ x\in B^m_\imath(n)\ {\rm such \ that}\ (1-\epsilon)\Gamma\\[1mm]
&\hskip-5mm \leq n^{-d}|\Lambda_n(x)|\leq (1+\epsilon)\Gamma\ \ {\rm does\  not\ hold}\ea\)\leq\exp(-c'_5 n^{1/2})$$ for large $n$. Where $c'_5$ be some positive constant less than $c_5=c_5(\a',\beta',\epsilon/2)$ given in (\ref{l2.01}).
Hence $$\P\(\ba{rl}& \hskip-5mm{\rm there \ exists}\ u\in V_n\ {\rm such\ that}\ (1-\epsilon)\Gamma\\[1mm]&\hskip-5mm \leq n^{-d}|\Lambda_n(u)|\leq (1+\epsilon)\Gamma\ {\rm does\ not\ hold}\ea \)\leq m^d\cdot\exp(-c'_5 n^{1/2})\rightarrow 0$$ as $n\rightarrow \infty$. \QED

For any positive integer $k$, and for any $\jmath=(j_1,j_2,\ldots,j_d)\in\{0,1,\ldots,k-1\}^d$, let $B^k_\jmath(n)$ be the box with side length $n/k$ centered at $x_\jmath:=\frac{n}{k}\jmath+(\frac n{2k},\frac n{2k},\ldots,\frac n{2k})$ and we call it a $k$-box. For any $\jmath=(j_1,j_2,\ldots,j_d)\in\{0,1,\ldots,k-1\}^d$,
let $\bar B^k_\jmath(n)=\frac n{k}\jmath+[\frac n{4k},\frac {3n}{4k}]^d$. Obviously, $B^k_\jmath(n)$ and $\bar B^k_\jmath(n)$ are all centered at $x_\jmath:=\frac nk\jmath+(\frac{n}{2k},\frac{n}{2k},\ldots,\frac{n}{2k})$ and can be rewritten as $B^{T_n}_{x_\jmath}(\frac n{2k})$ and $B^{T_n}_{x_\jmath}(\frac n{4k})$.
\bl\label{l3} Suppose $r>r_c$, and $0<\epsilon<1/2$. Let $k=k(n)=\lfloor{\hbar n}/{\ln^\psi n}\rfloor$ with $\hbar>0$ and $\psi\geq 1/(d-1)$. Then
\begin{description}
  \item[] (i) for small enough $\hbar>0$,
\be\label{l3.2}
\lim_{n\rightarrow\infty}\P\(|\bar B^k_\jmath(n)\cap V_n|\geq(1-\epsilon)\tilde\theta(r)\(n/2k\)^{d},\ \forall \ \jmath\ \rm{with}\ x_\jmath\in B_\varrho(n)\)=1,
\ee
where $B_\varrho(n)=[\varrho n,(1-\varrho)n]^d$ and $0<\varrho<1/2$ is a given constant;
  \item[] (ii) there exists $\delta=\delta(\epsilon)\in (0,1)$ such that
\be\label{l3.3}\lim_{n\rightarrow\infty}\P\(\left|\left\{\jmath: |B^k_\jmath(n)\cap V_n|\leq (1+2\epsilon)\tilde\theta(r)(n/k)^d\right\}\right|\geq \delta k^d\)=1;  \ {\rm and} \ee
 \item[] (iii) for any $\hbar>0$,
 \be\label{l3.4}\lim_{n\rightarrow\infty}\P\(|E(B^k_\jmath (n))|\leq 4(\hbar^{-1}\ln^\psi n)^{2d}, \ \forall\ \jmath\)=1,
\ee where $E(B^k_\jmath (n))$ is the set of edges of $G_n$ with both endpoints in $B^k_{\jmath}(n)$.
\end{description}
\el

{\it Proof.} {\it (i)} For any $\jmath$, denote by $\bar E^k_\jmath(n)$ the event related to box $\bar B^k_\jmath(n)$ as stated in Propositions~\ref{p1}. Since $n/k\geq \ln^\psi n/\hbar$, then, by Proposition~\ref{p1},
there exist $c'_1>0$ and $n'_0$ such that
\be\label{l3.01}
\P\(\bar E^k_\jmath(n)\)\geq 1-\exp\(-c'_1\hbar^{1-d}\ln^{(d-1)\psi} n\)\ee for all $n\geq n'_0$. Since $\hbar$ is small enough, then
\be\label{l3.001}\ba{rl}&\P\(\(\bigcap_\jmath \bar E^k_\jmath(n)\)^c\)\leq k^d\cdot\exp\(-c'_1\hbar^{1-d}\ln^{(d-1)\psi} n\)\\[2mm]&\leq \hbar^d \ln^{-d\psi} n \cdot n^d\cdot n^{-[c'_1\hbar^{1-d}\ln^{(d-1)\psi-1} n]}\\[2mm]&\rightarrow 0,\ \ \ \ {\rm as}\ n\rightarrow\infty. \ea \ee
Denote by $C_b(\bar B^k_\jmath(n))$ the unique largest cluster in $\bar B^k_\jmath(n)$ as stated in event $\bar E^k_\jmath(n)$, then (\ref{l3.001}) implies
\be\label{l3.01'}\P\(|C_b(\bar B^k_\jmath(n))|\geq(1-\epsilon)\tilde\theta(r)\(n/2k\)^{d}, \ \forall\ \jmath\)\rightarrow 1,\ \ \ {\rm as}\ n\rightarrow\infty;\ee and
\be\label{l3.01'''}\P\(C_b(\bar B^k_\jmath(n))\subset C_\infty, \ \forall\ \jmath\)\rightarrow 1,\ \ \ {\rm as}\ n\rightarrow\infty.\ee
Together with Proposition~\ref{p2}, (\ref{l3.01'''}) implies that
\be\label{l3.01''}\P\(V_n \supseteq C_b(\bar B^k_\jmath(n))\cap \fP^T_n, \forall\ \jmath\ {\rm with}\ x_\jmath\in B_\varrho(n)\)\rightarrow 1,\ \ {\rm as}\ n\rightarrow\infty.\ee
Thus, (\ref{l3.2}) follows from (\ref{l3.01'}) and (\ref{l3.01''}).

{\it (ii)} Let $l_\jmath:=|B^k_\jmath(n)\cap V_n|$, and $L(\epsilon):=\left|\left\{\jmath:l_\jmath\leq (1+2\epsilon)\tilde\theta(r)(n/k)^{d}\right\}\right|$. Then by Lemma~\ref{l1},
$$\lim_{n\rightarrow\infty}\P\([k^d-L(\epsilon)]\cdot (1+2\epsilon)\tilde\theta(r)(n/k)^{d}\leq (1+\epsilon)\tilde\theta(r) n^d\)=1.$$
This is just $\lim_{n\rightarrow\infty}\P\(L(\epsilon)\geq \delta(\epsilon)k^d\)=1$ with $\delta(\epsilon)=\epsilon/(1+2\epsilon)$.

{\it (iii)} For any $\jmath$, by the large deviation inequality (\ref{l1.05}), the probability that the number of Poisson points in $B^k_\jmath(n)$ exceed $2|B^k_\jmath(n)|$ is less than $\exp(-\gamma(2)(\hbar^{-1}\ln^\psi n)^d)$. Then, (\ref{l3.4}) follows immediately from the fact that $|E(B^k_\jmath(n))|$ can not exceed the square of the number of Poisson points in $B^k_\jmath(n)$.\QED

Finally, for $\Delta(G_n)$, the maximum degree of $G_n$, we have
\bl\label{l4} Suppose $r>r_c$. Then, for any $l\geq 2$,
\be\label{l4.1}\lim_{n\rightarrow\infty}\P\(\Delta(G_n)\leq \frac{\ln n}l\)=1.\ee\el

{\it Proof.} Let $k=\lfloor n/lr\rfloor$. $k$-boxes $B^k_{\imath}(n)$ and $B^k_{\jmath}(n)$ are called {\it adjacent}, if $|i_s-j_s|\leq 1$, or $=k-1$ for all $1\leq s\leq d$. Let $\hat B^k_{\imath}(n)$ denote the large box which is the union of $B^k_{\imath}(n)$ and all its neighbors in $T^d_n$. For any $x\in T^d_n$, let $Ba(x,lr)$ be the ball with radius $lr$ and centered at $x$.

If for some $u\in V_n$, $d_{G_n}(u)\geq \ln n/l$, then there exists $\imath\in\{0,1,2,\dots,k-1\}^d$ such that $u\in B^k_{\imath}(n)$, $Ba(u,lr)\subset \hat B^k_{\imath}(n)$ and
$$N(r):=|\hat B^k_{\imath}(n)\cap \fP^T_n|\geq \frac{\ln n}l.$$ By the large deviation inequality (\ref{l1.05}),
$$\ba{rl}&\P\({\rm there\ exists\ } u\in V_n{\rm\ such\ that}\ d_{G_n}(u)\geq \Df{\ln n}l\)\\[3mm]
&\Dp\leq k^d\P\(N(r)\geq \frac{\ln n}l\)\leq \(\frac n{lr}\)^d\cdot n^{-{c(l,r)\ln\ln n}}\\[3mm]
&\rightarrow 0, \ {\rm as\ } n\rightarrow\infty,\ea$$ where $c(l,r)>0$ is a constant depends on $l,r$.  \QED


\section{Geometry of $\fG_n$}
\renewcommand{\theequation}{3.\arabic{equation}}
\setcounter{equation}{0}

Before we give proofs to our main results, in this section, we shall study the geometry of $\fG_n$ in advance.

For any $r>r_c$, $\epsilon>0$ small enough (to choose $\epsilon$, see Remark~\ref{r1}), and for any $n$ large enough, we choose an arbitrary realization of the random graph $G_n=(V_n, E_n)$, still denote it by $G_n=(V_n,E_n)$, such that $G_n$ make all large probability events stated in Proposition~\ref{p1} and Lemmas~\ref{l1}-\ref{l4} occur. More precisely, with the notations given in the statements of Proposition~\ref{p1} and Lemmas~\ref{l1}-\ref{l4}, we assume the realization $G_n$, and the corresponding realizations of Poisson points $\fP^T_n$, $\fP_n$, satisfies the following conditions.
{\it
\begin{description}
  \item[\hskip 5mm $A_1:$] for $\fP_n$, $C_b(B(n))$ is crossing for $B(n)$ and is a part of the unbounded cluster $C_\infty$;
  \item[\hskip 5mm $A_2:$] for $V_n$, the vertex set of $G_n$, $(1-\epsilon)\tilde\theta(r)\leq n^{-d}|V_n|\leq (1+\epsilon)\tilde\theta(r)$;
  \item[\hskip 5mm $A_3:$] $(1-\epsilon)\Gamma\leq n^{-d}|\Lambda_n(u)|\leq (1+\epsilon)\Gamma$ for all $u\in V_n$;
  \item[\hskip 5mm $A_4:$] for some small enough $\hbar>0$ and $\varrho>0$, $|\bar B^k_\jmath(n)\cap V_n|$$\geq(1-\epsilon)\tilde\theta(r)\(n/2k\)^{d}$, for all $\jmath$ with  $x_\jmath\in B_\varrho(n)$, where $k=k(n)=\lfloor{\hbar n}/{\ln^\psi n}\rfloor$ with $\psi\geq 1/(d-1)$;
  \item[\hskip 5mm $A_5:$] $\left|\left\{\jmath: |B^k_\jmath(n)\cap V_n|\leq (1+2\epsilon)\tilde\theta(r)(n/k)^d\right\}\right|\geq \delta k^d$, where $k$ is given in $A_4$ and $\delta=\delta(\epsilon)=\epsilon/{(1+2\epsilon)}$;
  \item[\hskip 5mm $A_6:$] $|E(B^k_\jmath (n))|\leq 4(\hbar^{-1}\ln^\psi n)^{2d}$ for all $\jmath$, where $k$ is given in $A_4$;

   \item[\hskip 5mm $A_7:$] for some large enough integer $l\geq 2$, $\Delta(G_n)\leq \ln n/l$.

\end{description}
}

 Suppose that $\fG_n=\fG(G_n)$ is defined by adding random long edges to $G_n$. Clearly, by Proposition~\ref{p1} and Lemmas~\ref{l1}-\ref{l4}, it suffices to prove the main results of the paper for $\fG_n=\fG_n(G_n)$. We declare here that, in the rest of the paper, we only deal with $\fG_n=\fG_n(G_n)$ instead of $\fG_n$ as originally defined in Section 1.

On the deterministic realization $G_n=(V_n,E_n)$, for any vertex sets $S,S'\subset V_n$, let $\Lambda(S,S')$ denote the set of {\it unordered} vertex pairs $\{u,v\}$ with $u\in S$, $v\in S'$ and $v\in \Lambda_n(u)$, let $N(S,S')=|\Lambda(S,S')|$. To any {\it unordered} vertex pair $\{u,v\}\in \Lambda(V_n,V_n)$, independently, we assign a random variable $Z_{u,v}(=Z_{v,u})$ satisfying $\P(Z_{u,v}=1)=1-\P(Z_{u,v}=0)=p_n$. For any vertex sets $S,S'\subset V_n$, define
$$\fL(S,S'):=\sum_{\{u,v\}\in \Lambda(S,S')}Z_{u,v}.$$
Let $\fL(S):=\fL(S,S)+\fL(S,S^c)$. Clearly, $\fL(S,S)$ and $\fL(S,S^c)$ are independent binomial random variables with parameters $(N(S,S),p_n)$ and
$(N(S,S^c),p_n)$ respectively.
Let $N(S)=N(S,S)+N(S,S^c)$, then $\fL(S)$ is the binomial random variable with parameter $(N(S),p_n)$.

Define $\V(S):=\sum_{u\in S}d_{\fG_n}(u)$. Then
\be\label{3.5}\ba{rl}\V(S)&=\Dp\sum_{u\in S}d_{G_n}(u)+2\fL(S,S)+\fL(S,S^c)\\[3mm]&=\Dp\sum_{u\in S}d_{G_n}(u)+2\fL(S)-\fL(S,S^c).\ea\ee

To bound the tail probabilities of binomial random variable in the present paper, we introduce the following large deviation inequality.

\bl{Suppose $Z\sim b(n,p)$. Then
\be\label{3.15}\P(Z\geq zn)\leq \exp(-I(z)n),\ z>p \ {\rm and }\ \ \P(Z\leq zn)\leq \exp(-I(z)n), \ z<p,\ee where $I(z)$ is the common rate function defined by \be\label{3.16} I(z):=z\ln\frac{zq}{(1-z)p}-\ln \frac q{1-z}, \ p\not=z\in (0,1),\ \ p+q=1.\ee} Especially for small $p$, (\ref{3.15}) can be rewritten as
\be\label{3.15'}\ba{rl}&\P(Z\geq zpn)\leq \exp(-\gamma(z)pn)\ {\it for}\ z>1, \ {\rm and }\\[3mm]&\P(Z\leq zpn)\leq \exp\(-\frac 12\gamma(z)pn\)\ {\it for}\ 0<z<1,\ea\ee with $\gamma(z)=z\ln z-z+1$, same defined as in (\ref{l1.05}).\el

{\it Proof.} (\ref{3.15}) follows from the proof of the classical Cram\'er's Theorem \cite{C}. (\ref{3.15'}) follows from (\ref{3.15}) by using the Taylor's expansion of $I(zp)$ for small $p$.\QED

First of all, we shall bound  $\Delta(\fG_n)$, the {\it maximum degree} of $\fG_n$ from above. Actually, we have the following lemma.

\bl\label{le3.1} For $\chi> 1$ and $\sigma>0$, or for $\chi=1$ and $\sigma$ is large enough, there exists some constant $M=M(\sigma)$ large enough, such that \be\label{3.1}\lim_{n\rightarrow\infty}\P(\Delta(\fG_n)\leq M\ln^{\chi} n)=1.\ee\el

{\it Proof.} For any $u\in V_n$, let $\fL(u):=\fL(\{u\})=d_{\fG_n}(u)-d_{G_n}(u)$. Then by conditions $A_2$ and $A_3$,
$$\P\(\max_{u\in V_n}\fL(u)\geq M_1 \ln^{\chi} n\)\leq(1+\epsilon)\tilde\theta(r)n^d\P\(\xi\geq M_1 \ln^{\chi} n\),$$ where $M_1=2(1+\epsilon)\Gamma\sigma$, $\xi$ is the Binomial random variable with parameter $((1+\epsilon)\Gamma n^d, p_n)$ and $\Gamma$ is the constant defined in Lemma~\ref{l2}. Using the large deviation inequality (\ref{3.15'}), we have $\P\(\xi\geq M_1 \ln^{\chi} n\)\leq n^{-\Omega(\sigma \ln^{\chi-1} n)}$. Then,
\be\label{3.3}\lim_{n\rightarrow\infty}\P\(\max_{u\in V_n}\fL(u)\leq M_1\ln^{\chi} n\)=1.\ee
The lemma follows from condition $A_7$, (\ref{3.3}) and the fact that $\Delta(\fG_n)\leq \Delta(G_n)+\Dp\max_{u\in V_n}\fL(u)$.\QED

For the lazy random walk $\{X_t:t\geq 0\}$ on $\fG_n$, let $Q(u,v):=\pi(u)P(u,v)$ and $Q(S,S^c):=\sum_{u\in S}\sum_{v\in S^c}Q(u,v)$. Define
$$h:=\min_{S:\pi(S)\leq 1/2}\frac{Q(S,S^c)}{\pi(S)}$$ to be the {\it conductance} of $\{X_t:t\geq 0\}$. Letting $\fE(S,S^c)$ be the number of edges between $S$ and $S^c$, we have \be\label{3.6}h=\frac 12\min_{S:\pi(S)\leq 1/2}\frac{\fE(S,S^c)}{\V(S)}.\ee Obviously, the conductance $h$ is ultimately determined by the geometry of $\fG_n$.

Another interesting quality on $\fG_n$ is the {\it edge isoperimetric constant} $\iota$ defined by
\be\label{3.7}\iota:=\min_{S:|S|\leq |V_n|/2}\frac{\fE(S,S^c)}{|S|}.\ee Note that the edge isoperimetric constant of a graph is also called the {\it Cheeger constant} in honor of the eigenvalue bound in differential geometry. In the rest of this section, we will try to give lower bounds to $h$ and $\iota$. Using these lower bounds, we then finish the proofs of our main results in the next section.


\bl\label{le3.2}Suppose $\chi\geq 1$. Then for small enough $a>0$, we have
\be\label{3.4}\lim_{n\rightarrow\infty}\P\(\bigcap_{S:|S|\geq (1-a)|V_n|}\left\{\pi(S)>\frac 12\right\}\)=1.\ee \el

{\it Proof.} First of all, for any $S\subset V_n$, by condition $A_7$, we have
\be\label{3.11'}\pi(S)=\Df{\sum_{u\in S}d_{G_n}(u)+2\fL(S,S)+\fL(S,S^c)}{\sum_{u\in V_n}d_{G_n}(u)+2\fL(V_n)}\geq\Df{2\fL(S,S)}{|V_n|\ln n/l+2\fL(S,S)+2\fL(S^c)}.\ee


For any given $S$ with $|S|\geq (1-a)|V_n|$, by definitions, $\fL(S,S)\sim b(N(S,S),p_n)$, $\fL(S^c)\sim b(N(S^c), p_n)$ are independent binomial random variables, and for small enough $a>0$,
\be\label{3.10}N(S^c)\leq [(1-a)a+a^2/2]\cdot |V_n|^2\leq\epsilon_1 A|V_n|^2,\ee  \be\label{3.11}N(S,S)\geq \left[\frac{(1-\epsilon)\Gamma}{2(1+\epsilon)\tilde\theta(r)}-\(1-a\)a-\frac{a^2}2\right]|V_n|^2\geq (1-\epsilon_1)A |V_n|^2,\ee where $A=\frac{(1-\epsilon)\Gamma}{2(1+\epsilon)\tilde\theta(r)}$ and $\epsilon_1>0$ is a given small constant.
 Using the inequality (\ref{3.15'}), we know that both $\P\(\fL(S,S)\leq (1-2\epsilon_1)A|V_n|^2p_n\)$ and $\P\(\fL(S^c)\geq 2\epsilon_1 A|V_n|^2p_n\)$ are less than $\exp\(-\Omega(n^d\ln^{\chi} n)\)$. Note that, by (\ref{3.11'}), the fact that $\chi\geq 1$, $\fL(S,S)\geq (1-2\epsilon_1)A|V_n|^2p_n$ and $\fL(S^c)\leq 2\epsilon_1|V_n|^2p_n$ imply that
$$\pi(S)\geq\Df{2l\fL(S,S)}{|V_n|\ln n+2l\fL(S,S)+4l\epsilon_1 A|V_n|^2p_n}\geq\Df{2l(1-2\epsilon_1)A|V_n|p_n}{\ln n+2lA|V_n|p_n}> \frac 12,$$ for large $n$, $l$ and small $\epsilon_1$ (here we only need $l$ large enough in case of $\chi=1$).
Then, we obtain
\be\label{3.13}\P\(\pi(S)>\frac 12\)\geq 1-\exp\(-\Omega(n^{d}\ln^{\chi}n)\).\ee
Now, let $M_a:=|\{S\subset V_n:|S|\geq (1-a)|V_n|\}|$. To finish the proof of the lemma, it remains to bound $M_a$ from above.

By Lemma 6.3.3 in \cite{D}, the number of $S\subset V_n$ with $|S|=s$ is
\be\label{3.25'}\b{|V_n|}{s}\leq \(\Df{|V_n|\cdot e}{s}\)^s=\exp\left\{s\left[\ln\(\frac{|V_n|}s\)+1\right]\right\}.\ee
Then
\be\label{3.18}M_a=\Dp\sum_{s=(1-a)|V_n|}^{|V_n|}\b{|V_n|}{s}\leq a|V_n|\cdot\b{|V_n|}{a|V_n|}\leq \exp \(O\(a[\ln{(1/a)}+1]n^d\)\).\ee
Combining (\ref{3.13}) and (\ref{3.18}), for small enough $a>0$, we obtain $$\P\(\bigcup_{S:|S|\geq (1-a)|V_n|}\left\{\pi(S)\leq\frac 12\right\}\)\leq M_a\cdot\exp\(-\Omega(n^{d}\ln^{\chi}n)\)\rightarrow 0$$ as $n\rightarrow\infty$.\QED

Let \be\label{3.8}\hat \iota=\min_{S:|S|\leq(1-a)|V_n|}\frac{\fE(S,S^c)}{\V(S)}.\ee Then, by Lemma~\ref{le3.2},
\be\label{3.88}\P\(h\geq \hat \iota/2\)\rightarrow 1,\ \ {\rm as}\ n\rightarrow \infty.\ee So, to bound $h$ from below, it suffices to bound $\hat \iota$ from bellow. In fact, we have

\bp\label{p3.3}Suppose $0<\a<\beta<1/2$ with $(2\beta)^d-(2\a)^d>1/2$, then, for $\chi>1 $ and $\sigma>0$, or for $\chi=1$ and $\sigma$ is large enough, there exists $C_4>0$ such that
\be\lim_{n\rightarrow\infty}\P\(\hat \iota\geq {C_4}{\ln^{-1} n}\)=1.\ee\ep

{\it Proof.} The proof of this proposition is the main part of our proofs. In fact, we will develop a more complicated version of the approach proposed by Durrett in \cite[Theorem 6.6.1]{D}. Note that in \cite{D}, the mixing times of random walks on several small worlds were studied.

Let
$$\ba{rl}&\hB_1:=\{S\subset V_n: |S|\leq a|V_n|\},\ {\rm and}\\[3mm]
&\hB_2:=\{S\subset V_n:a|V_n|<|S|\leq (1-a)|V_n|\},\ea$$ where $a>0$ is given in Lemma~\ref{le3.2} and we also also assume that $(2\beta)^d-(2\a)^d>1/2+a$. The proposition follows from the
following Lemmas~\ref{le3.4} and \ref{le3.6}.\QED

\bl\label{le3.4} Suppose $0<\a<\beta<1/2$. For $\chi>1$ and $\sigma>0$, or for $\chi=1$ and $\sigma$ is large enough, there exists $C_5>0$ such that
\be\label{3.17}\lim_{n\rightarrow\infty}\P\left\{\Df{\fE (S,S^c)}{\V(S)}\geq\Df{C_5}{\ln n},\ \forall\ S\in\hB_1\right\}=1.\ee\el

{\it Proof.} For any $S\subset V_n$, let's consider the random variables $\fL(S,S^c)$, $\fL(S,S)$ and $\fL(S)$. Recall that $\fL(S,S^c)\sim b(N(S,S^c),p_n)$, $\fL(S,S)\sim b(N(S,S),p_n)$ and $\fL(S)\sim b(N(S),p_n)$.

First of all, we have
\be\label{3.24}\ba{rl}&N(S)\leq|S|\cdot\max_{u\in S}|\Lambda_n(u)|\leq (1+\epsilon)\Gamma n^d|S|=:N_1(S),\\[3mm]
&N(S)\geq\frac 12|S|\cdot\min_{u\in S}|\Lambda_n(u)|\geq \frac12(1-\epsilon)\Gamma n^d|S|=:N_2(S), \ea\ee note that this indicates that $N(S)=\Theta(n^d|S|)$.
Then, by definition and (\ref{3.24})
$$\ba{rl}N(S,S^c)&\geq \Dp\sum_{u\in S}\(\sum_{v\in \Lambda_n(u)}1-|S|\)\geq |S|\((1-\epsilon)\Gamma n^d-a|V_n|\)\geq (1-2\epsilon)\Gamma n^d|S|\\[8mm]&=\Df{1-2\epsilon}{1+\epsilon}N_1(S)\geq \Df{1-2\epsilon}{1+\epsilon}N(S). \ea$$ Hence $$N(S,S)=N(S)-N(S,S^c)\leq \(1-\frac{1-2\epsilon}{1+\epsilon}\)N(S)=\frac{3\epsilon}{1+\epsilon}N(S).$$

Let $W_1(S)\sim b(3\epsilon N(S)/(1+\epsilon),p_n)$, $W_2(S)\sim b((1-2\epsilon)N(S)/(1+\epsilon),p_n)$. Suppose that $W_1(S)$ and $W_2(S)$ are independent. Then,
$$\ba{rl}\P\(\fL(S,S^c)\geq \Df{\fL(S)}2\)&=\P(\fL(S,S^c)\geq \fL(S,S))\geq \P(W_2(S)\geq W_1(S))\\[4mm]&\geq \P\(W_2(S)> \Df {N(S)p_n}2> W_1(S)\).\ea$$
Using the inequality (\ref{3.15'}) and the fact that $N(S)\geq N_2(S)=\frac 12 (1-\epsilon)\Gamma n^d |S|$, we obtain that both $\P\(W_2(S)\leq N(S)\cdot p_n/2\)$ and $\P\(W_1(S)\geq N(S)\cdot p_n/2\)$ are less than $\exp\(-\Omega\(\sigma |S|\ln^{\chi} n\)\)$. Then
\be\label{3.28}\P\(\fL(S,S^c)\geq \frac{\fL(S)}2\)\geq 1-\exp\(-\Omega\(\sigma|S|\ln^{\chi} n\)\).\ee

Let $Z_2(S)\sim b(N_2(S),p_n)$, then
$$\P (\fL(S)\leq |S|)\leq \P(Z_2(S)\leq |S|)=\P\(Z_2(S)\leq \frac{2\ln^{-\chi} n}{\sigma(1-\epsilon)\Gamma}N_2(S)p_n\).$$
Using the inequality (\ref{3.15'}) again, we obtain
\be\label{3.30}\P\(\fL(S)\leq |S|\)\leq\exp\(-\Omega\(\sigma|S|\ln^{\chi} n\)\).\ee
Note that $\fL(S,S^c)\geq {\fL(S)}/2$, $\fL(S)\geq |S|$ and condition $A_7$ imply

\be\label{3.32}\Df{e(S,S^c)}{\V(S)}\geq \Df{\fL(S,S^c)}{|S|\Delta(G_n)+2\fL(S)}
\geq\Df{\fL(S)/2}{\fL(S)\ln n/2+2\fL(S)}\geq\Df{C_5}{\ln n}.\ee



Using (\ref{3.25'}), (\ref{3.28}) and (\ref{3.30}), we obtain \be\label{3.31}\ba{rl}&\P\(\Dp\bigcup_{S\in \hB_1}\{\fL(S,S^c)\geq {\fL(S)}/2\ {\rm and\ }\fL(S)\geq |S|\}^c\)\\[4mm]&\leq\Dp\sum_{s=1}^{a|V_n|}\b{|V_n|}{s}\Dp\exp\(-\Omega(\sigma s\ln^{\chi} n)\)\\[6mm]
&\leq \Dp n^{d}\exp\(-[\Omega(\sigma \ln^\chi n)-O(\ln n)]\)\rightarrow 0,\ {\rm as} \ n\rightarrow\infty.\ea\ee
The lemma follows immediately from (\ref{3.32}) and (\ref{3.31}).\QED

\bl\label{le3.6} Suppose $0<\a<\beta<1/2$ with $(2\beta)^d-(2\a)^d>1/2$, then for $\chi>0$ and $\sigma>0$, or for $\chi=0$ and $\sigma$ is large enough, there exist $C_6>0$ such that
\be\label{3.22}\lim_{n\rightarrow\infty}\P\left\{\Df{\fE (S,S^c)}{\V(S)}\geq\Df{C_6}{\ln n},\ \forall\ S\in\hB_2\right\}=1.\ee\el

{\it Proof.} The proof of this lemma is very similar to the proof of Lemma~\ref{le3.4}. But just in this step, we have to use the condition of $\a$ and $\beta$: $(2\beta)^d-(2\a)^d>1/2.$

%

Recall that $a (>0)$ is given in Lemma~\ref{le3.2} and is chosen small enough such that $(2\beta)^d-(2\a)^d>1/2+a$. Now we choose $\epsilon>0$ small enough such that
$$ \Upsilon:=\Df{1-\epsilon}{(1+\epsilon)\tilde\theta(r)}\Gamma=\Df{(1-\epsilon)[(2\beta)^d-(2\a)^d]}{(1+\epsilon)}\geq \frac 12+a,$$ and then $1-1/(2\Upsilon)\geq a/\Upsilon$.

For any $S\in \hB_2$, if $a|V_n|\leq |S|\leq |V_n|/2$, then
$$\ba{rl}N(S,S^c)&=\Dp\sum_{u\in S}\sum_{v\in S^c\cap\Lambda_n(u)}1\geq |S|\(\min_{u\in S}|\Lambda_n(u)|-|S|\)\geq |S|\(\Upsilon|V_n|-|S|\)\\[4mm]
&\geq a|V_n|\(\Upsilon|V_n|-a|V_n|\)\geq a\(\Df{1-2\epsilon}{(1+\epsilon)\tilde\theta(r)}\Gamma\)|V_n|^2,\ea$$
where the third inequality comes from the fact that the function $g(x)=x(1-x)$ in interval $[a/\Upsilon,1/(2\Upsilon)]$ takes its minimum at $x=a/\Upsilon$. On the other hand, if $a|V_n|\leq |S^c|\leq |V_n|/2$, then
$$\ba{rl}N(S,S^c)&=N(S^c,S)=\Dp\sum_{u\in S^c}\sum_{v\in S\cap\Lambda_n(u)}1\geq |S^c|\(\min_{u\in S^c}|\Lambda_n(u)|-|S^c|\)\\[3mm]
&\geq a|V_n|\(\Upsilon|V_n|-a|V_n|\)\geq a\(\Df{1-2\epsilon}{(1+\epsilon)\tilde\theta(r)}\Gamma\)|V_n|^2.\ea$$
So, by condition $A_2$, for any $S\in \hB_2$,
$$\ba{rl}N_1(S)&=(1+\epsilon)\Gamma n^d |S|\leq (1+\epsilon)\Gamma n^d (1-a)|V_n|\leq \Df{(1+\epsilon)\Gamma(1-a)}{(1-\epsilon)\tilde\theta(r)}|V_n|^2\\[3mm]
&\leq \Df{(1+\epsilon)\Gamma }{(1-\epsilon)\tilde\theta(r)}|V_n|^2\leq\frac1a\cdot\Df {(1+\epsilon)^2}{(1-\epsilon)(1-2\epsilon)}N(S,S^c).\ea$$
Let $f(\epsilon)={(1-\epsilon)(1-2\epsilon)}/{(1+\epsilon)^2}$, then \be\label{3.35} N(S,S^c)\geq f(\epsilon)a N_1(S)\geq f(\epsilon)a N(S).\ee

Now, by the large deviation inequality (\ref{3.15'}), we have
$$\ba{rl}&\P\(\fL(S,S^c)\leq \frac 12f(\epsilon)aN(S)p_n\)\\[2mm]&\leq\P_\epsilon\(\fL(S,S^c)\leq \frac 12N(S,S^c)p_n\)\leq \exp\(-\frac12\gamma\(\frac 12\)N(S,S^c)p_n\)\\[2mm]&\leq\exp\(-\frac12\gamma\(\frac 12\)f(\epsilon)aN(S)p_n\)\leq\exp\(-\frac12\gamma\(\frac 12\)f(\epsilon)aN_2(S)p_n\)\\[2mm]
&\leq\exp(-\Omega(\sigma |S|\ln^{\chi} n))\ea$$ and
$$\ba{rl}\P\(\fL(S)\geq 2N(S)p_n\)&\leq \exp\(-\gamma(2)N(S)p_n\)\leq  \exp\(-\gamma(2)N_2(S)p_n\)\\[2mm]
&\leq\exp(-\Omega(\sigma |S|\ln^{\chi} n)).\ea$$
If $\fL(S,S^c)\geq \frac 12f(\epsilon)aN(S)p_n$ and $\fL(S)\leq 2N(S)p_n$, then
$\fL(S,S^c)\geq\frac 14f(\epsilon)a \fL(S)$. So
\be\label{3.36}\ba{rl}&\P\(\fL(S,S^c)\geq f(\epsilon) a\fL(S)/4\)\\[2mm]&\geq \P(\fL(S,S^c)\geq \frac 12f(\epsilon)aN(S)p_n \ {\rm and}\ \fL(S)\leq 2N(S)p_n)\\[2mm]
&\geq 1-\exp\(-\Omega\(\sigma|S|\ln^{\chi}n\)\).\ea\ee
Now $\fL(S,S^c)\geq f(\epsilon) a\fL(S)/4$, $\fL(S)\geq |S|$ and condition $A_7$ imply

\be\label{3.37}\Df{e(S,S^c)}{\V(S)}\geq \Df{\fL(S,S^c)}{|S|\Delta(G_n)+2\fL(S)}
\geq\Df{f(\epsilon) a\fL(S)/4}{\fL(S)\ln n/2+2\fL(S)}\geq\Df{C_7}{\ln n}.\ee


Note that (\ref{3.30}) also holds for $S\in \hB_2$. Using (\ref{3.25'}), (\ref{3.30}) and (\ref{3.36}), we obtain \be\label{3.38}\ba{rl}&\P\(\Dp\bigcup_{S\in \hB_2}\{\fL(S,S^c)\geq f(\epsilon) a\fL(S)/4\ {\rm and\ }\fL(S)\geq |S|\}^c\)\\[6mm]&\leq\Dp\sum_{s=a|V_n|}^{(1-a)|V_n|}\b{|V_n|}{s}\Dp\exp\(-\Omega(\sigma s\ln^{\chi} n)\)\\[6mm]&\leq \Dp\sum_{s=a|V_n|}^{(1-a)|V_n|}\exp\(-[\Omega(\sigma \ln^\chi n)-(\ln (1/a)+1)]s\)\rightarrow 0,\ {\rm as} \ n\rightarrow\infty.\ea\ee
The Lemma now follows from (\ref{3.37}) and (\ref{3.38}).\QED

Similar to Proposition~\ref{p3.3}, we can obtain the following lower bound for $\iota$, the edge isoperimetric constant of $\fG_n$.

\bp{\label{p3.8}} Suppose $0<\a<\beta<1/2$ with $(2\beta)^d-(2\a)^d>1/2$, then, for $\chi> 1$ and $\sigma>0$, or for $\chi=1$ and $\sigma$ is large enough, there exists $C_7>0$ such that
\be\lim_{n\rightarrow\infty}\P\(\iota\geq {C_7}\)=1.\ee\ep

{\it Proof.} Let
$$\ba{rl}&\hB'_1:=\{S\subset V_n: |S|\leq a|V_n|\},\ {\rm and}\\[3mm]&\hB'_2:=\{S\subset V_n:a|V_n|<|S|\leq |V_n|/2\},\ea$$ where $a>0$ is small enough such that $(2\beta)^d-(2\a)^d>1/2+a$. The proposition follows immediately from the inequality
$$\Df {e(S,S^c)}{|S|}\geq\Df{\fL(S,S^c) }{|S|}$$ and the estimates given in Lemmas \ref{le3.4} and \ref{le3.6}.\QED

\br\label{r1} We determine the $\epsilon$ given in the beginning of Section 3 as follows. Firstly, as in the proof of Lemma~\ref{le3.2}, we choose $l$ (in $A_7$) large enough such that we can determine a small enough $a>0$ uniformly for $\epsilon\in(0,1/2)$. Secondly, for given $0<\a<\beta<1/2$ with $(2\beta)^d-(2\a)^d>1/2$, reset $a$ such that we also have $(2\beta)^d-(2\a)^d>1/2+a$. Finally, we choose $\epsilon$ small enough such that $\Upsilon:={(1-\epsilon)[(2\beta)^d-(2\a)^d]}/{(1+\epsilon)}\geq \frac 12+a$. \er

\section{Proofs of Theorems~\ref{th1} and \ref{th2}}
\renewcommand{\theequation}{4.\arabic{equation}}
\setcounter{equation}{0}

Recall that $G_n=(V_n,E_n)$ is the realization of the random graph given in the beginning of Section 3, and $\fG_n=\fG_n(G_n)$ is the resulting random graph by adding random long edges to $G_n$. In this section, we will prove our mail results for $\fG_n=\fG_n(G_n)$.

{\it Proof of (ii) of Theorem~\ref{th1}.} Noticing that ${\rm diam}(\fG_n)$ is non-increasing in $\chi$, we obtain item (ii) of Theorem~\ref{th1} by Lemma~\ref{le3.1}, Proposition~\ref{p3.8} and the following Lemma~\ref{le4.1}.\QED

\bl\label{le4.1} For any connected graph $G=(V,E)$, let $\Delta(G)$ denote its maximum degree, $\iota(G)$ denote its edge isoperimetric constant and let ${\rm diam(G)}$ denote its diameter. Then \be\label{4.1}{\rm diam(G)}\leq \frac{4\Delta(G)}{\iota(G)}\ln |V|.\ee\el

{\it Proof.} This is a well known result in algebraic graph theory, for a detailed proof, one may refer to \cite{AM,BZ}.\QED

{\it Proofs of (i) of Theorem~\ref{th1} and (i) of Theorem~\ref{th2}.} In this part of proof, let's recall that our $G_n$ satisfies the conditions $A_4$-$A_6$.

Fix an integer $m_0$ such that $0<1/m_0<\a$, let's consider the family of $m_0$-boxes: $B^{m_0}_\imath(n)$, ${\imath}=(i_1,i_2,\ldots,i_d)\in\{0,1,2,\ldots, m_0-1\}^d$. Let $k=k(n)=\lfloor{\hbar n}/{\ln^\psi n}\rfloor$ with $\psi=\frac{1-\chi}d\geq \frac 1{d-1}$, where $\hbar>0$ is given in condition $A_4$. A $k$-box $B^k_\jmath(n)$, $\jmath\in\{0,1,\ldots,k-1\}^d$, is called {\it good}, if its center is in $B_\varrho(n)$ and $l_\jmath=|B^k_\jmath(n)\cap V_n|\leq (1+2\epsilon)\tilde\theta(r)(n/k)^{d}$, where $\varrho>0$ is given in condition $A_4$. A $m_0$-box $B^{m_0}_\imath(n)$ is called {\it good}, if it intersects more than $\frac 12\delta(k/m_0)^d$ good $k$-boxes, where $\delta=\delta(\epsilon)$ is given in condition $A_5$.

On one hand, by condition $A_5$,  for $G_n$, there exists a $m_0$-box which intersects more than $\delta(k/m_0)^d$ $k$-boxes, and each of these $k$-boxes contains at most $(1+2\epsilon)\tilde\theta(r)(n/k)^{d}$ vertices in $V_n$. On the other hand, by condition $A_4$, $\varrho>0$ is small enough, then the ratio of $k$-box with its center {\it not} in $B_\varrho(n)$ can be arbitrary small. So, the above existed $m_0$-box will intersects more than $\frac 12 \delta(k/m_0)^d$ good $k$-boxes, and it is really a {\it good} $m_0$-box.
Namely, for $G_n$, a good $m_0$-box exists.

 Now, suppose that $B^{m_0}_{\imath_0}(n)$ is a good $m_0$-box. Let $J(n)$ denote the set of index $\jmath$ of good $k$-box $B^k_\jmath(n)$ which intersects $B^{m_0}_{\imath_0}(n)$.  For any $\jmath\in J(n)$, let $L_\jmath$ be the {\it random} number of {\it long edge} (in $\fG_n$) with one of its endpoints in $B^k_\jmath(n)$. Because the side length of an $m_0$-box is $n/m_0$, and the side length of an $k$-box is $\Omega(\ln^{\psi}n)$, then for any $x\in B^k_{\jmath_1}(n)$, $y\in B^k_{\jmath_2}(n)$, $\jmath_1,\jmath_2\in J(n)$, $\jmath_1\not=\jmath_2$, one has
 $$d^T_\infty(x,y)\leq n/m_0+2\Omega(\ln^{\psi}n)<\a n, \ \rm{for}\ \ \rm{large} \ n.$$ Hence, the long edges counted in $L_{\jmath_1}$ differ from the long edges counted in $L_{\jmath_2}$ and random variables $\{L_\jmath:\jmath\in J(n)\}$ are independent.

 For any $\jmath\in J(n)$, let $\lambda_\jmath:=\Dp\max_{u\in B^k_\jmath(n)\cap V_n}|\Lambda_n(u)|$, then by condition $A_3$, $\lambda_\jmath\leq(1+\epsilon)\Gamma n^d$. Hence, \be\label{4.2}\ba{rl}\P(L_\jmath=0)&\geq(1-p_n)^{l_\jmath\cdot\lambda_\jmath}\geq(1-p_n)^{(1+2\epsilon)\tilde\theta(r)(n/k)^d\cdot(1+\epsilon)\Gamma n^d}\\[3mm] &\geq (1-p_n)^{\eta n^d\ln^{d\psi} n}\geq \exp\({-2\sigma\eta\ln^{d\psi+\chi}n}\)\\[3mm]
 &=n^{-2\sigma\eta},\ea \ee where $\eta=(1+2\epsilon)(1+\epsilon)\tilde\theta(r)\Gamma 2^d\hbar^{-d}$. Note that in the last inequality, we used the fact that $1-x>e^{-2x}$ for small $x>0$. Then, by condition $A_5$ and the independence proved above,
\be\label{4.3}\ba{rl}&\P\(\bigcap_{\jmath\in J(n)}\{L_\jmath\geq 1\}\)=\prod_{\jmath\in J(n)}\P(L_\jmath\geq 1)\leq (1-n^{-2\sigma\eta})^{\frac 12\delta(k/m_0)^d}\\[4mm]
 &\leq\exp\{-O(n^{d-2\sigma\eta}\ln^{-d\psi} n)\},\ {\rm for\ large}\ n.\ea \ee

By taking $\sigma>0$ small enough, we always have $2\sigma\eta<d$. Then, by (\ref{4.3}), \be\label{4.8}\lim_{n\rightarrow\infty}\P\(\bigcup_{\jmath\in J(n)}\{L_\jmath=0\}\)=1.\ee
Let $B^k_{\jmath_0}(n)$ be a good $k$-box with $L_{\jmath_0}=0$. By conditions $A_4$ and $A_6$, we also have $\bar B^k_{\jmath_0}(n)\cap V_n \not=\phi$ and $|E^k_{\jmath_0}(n)|\leq 4(\hbar^{-1}\ln^\psi n)^{2d}$. Suppose $u_0\in \bar B^k_{\jmath_0}(n)\cap V_n$, then for any $u\in [B^k_{\jmath_0}(n)]^c\cap V_n$, ${\cal D}_{\fG_n}(u_0,u)\geq \ln^\psi n/8r\hbar$, this first finishes the proof of (i) of Theorem~\ref{th1}.

Now, let's consider the lazy random walk on $\fG_n$. Clearly, the random walk started from $u_0$ can not escape from the box $B^k_{\jmath_0}(n)$ in time $C_2\ln^\psi n$, where $C_2=1/8r\hbar$. Then, for any $t\leq C_2\ln^\psi n$, $P^t(u_0, u)=0$ for all $u\in [B^k_{\jmath_0}(n)]^c\cap V_n$. So, by condition $A_6$,
 $$\ba{ll} ||&\hskip-4mm P^t(u_0,\cdot)-\pi||_{TV}\geq \frac 12\pi\([B^k_{\jmath_0}(n)]^c\cap V_n\)\geq \frac 12-\frac{|E(B^k_{\jmath_0}(n))|}{2|E_n|}\\[2mm]&\hskip-3mm\geq \frac 12-\frac{8r(\hbar^{-1}\ln^\psi n)^{2d}}{2n}>\frac 1 e, \ea$$ for large $n$. Note that in the third inequality, we have used the fact that $|E_n|\geq n/2r$, which follows from condition $A_1$. By definition of mixing time, we have $T_{\rm mix}\geq C_2\ln^\psi n$ and finish the proof of (i) of Theorem~\ref{th2}.\QED


 %

 {\it Proof of (ii) of Theorem~\ref{th2}.} For our lazy random walk $\{X_t:t\geq 0\}$ on $\fG_n$, matrix theory tell us that the transition kernel (P(u,v)) has nonnegative real eigenvalues $$1=\lambda_0\geq\lambda_1\geq\lambda_2\geq\ldots\geq\lambda_{|V_n|-1}\geq 0.$$ Note that $1-\lambda_1$ is called the spectral gap of $(P(u,v))$. Let $\pi_{\rm min}=\min_{u\in V_n}\pi(u)$.

 As a standard relation, it can be found in \cite[Theorem 12.5]{LPW} that
 \be\label{4.6}T_{\rm mix}\leq\ln\(\frac e{\pi_{\rm min}}\)\frac 1{1-\lambda_1}.\ee

The spectral gap $1-\lambda_1$ can be bounded from above and below by the conductance $h$ in the following way (see \cite[Theorem 6.2.1]{D}),
\be\label{4.5}\frac{h^2}2\leq 1-\lambda_1\leq 2h.\ee
On the other hand, it follows from Lemma~\ref{le3.1} that
\be\label{4.7}\P\(\Df 1{\pi_{\rm min}}\leq M|V_n|\ln^{\chi} n\)\geq \P\(\Delta(\fG_n)\leq M\ln^{\chi} n\)\rightarrow 1, \ {\rm as\ }n\rightarrow\infty.\ee
Thus, the desired result follows from (\ref{3.88}), (\ref{4.6})-(\ref{4.7}) and Proposition~\ref{p3.3}.\QED

\section*{Acknowledgements}
The author much thank one of the
anonymous referees for pointing out an important correction
in the proof of Lemma~\ref{le3.4} in an earlier version,
as well as for subsequent
discussions on our fixing of them.
The author also thank the
anonymous referees for other valuable comments and suggestions.

\vskip10mm

\begin{minipage}{6.5cm}

\noindent School of Mathematical Sciences, Capital
Normal University, Beijing, 100048, China. Email:
\ \texttt{wuxy@cnu.edu.cn}

\end{minipage}

\end{document}